\documentclass[12pt,twoside,letterpaper]{article}

\usepackage{fancyhdr}
\usepackage{natbib}
\usepackage[latin1]{inputenc}
\usepackage[T1]{fontenc}
\usepackage{amsmath}
\usepackage{amssymb}
\usepackage{graphicx}
\usepackage{amsthm,color}

\newfam\bbfam
\font\tenbb   =msbm10 \textfont\bbfam=\tenbb
\font\sevenbb =msbm7  \scriptfont\bbfam=\sevenbb
\font\fivebb  =msbm5  \scriptscriptfont\bbfam=\fivebb

 \newcommand{\eps}{\varepsilon}
 \newcommand{\To}{\longrightarrow}
 \newcommand{\B}{\mathcal{B}}

 \newcommand{\E}{\mathbb{E}}
 \newcommand{\X}{\mathcal{X}}
 \newcommand{\F}{\mathcal{F}}

 \newcommand{\Real}{\mathbb{R}}

 \newcommand{\abs}[1]{\left\vert#1\right\vert}

 \newcommand{\norm}[1]{\left\Vert#1\right\Vert}
 
\newcommand{\nnorm}[1]{\left\vert\!\left\vert\!\left\vert#1\right\vert\!\right\vert\!\right\vert}

\newtheorem{theorem}{Theorem}[section]

\newtheorem{algo}{Algorithm}[section]
\newtheorem{lemma}{Lemma}[section]

\newcounter{saveeqn}

\newlength{\defaultheadheight}
\newlength{\defaultheadsep}
\headheight \defaultheadheight
\begin{document}

\pagenumbering{arabic}
\pagestyle{fancy}
\fancyhf{} 
\renewcommand{\footrulewidth}{0pt}
\renewcommand{\headrulewidth}{0pt}

\fancyhead[CO]{\scriptsize Resampling from the past}
\fancyhead[RO]{\thepage}
\fancyhead[CE]{\scriptsize Y. F. Atchad\'e}
\fancyhead[LE]{\thepage}

\thispagestyle{empty}

\begin{center}
{\Large {\bf Resampling from the past to improve on MCMC algorithms\footnote{This work is funded in part by NSERC Canada}}} \\[12pt]

{\footnotesize {\em Yves F. Atchad\'e\footnote{Department of
Mathematics and Statistics, University of Ottawa, email:
yatchade@uottawa.ca}}\\[5pt]
(May 2006)} \\[10pt]
\end{center}

\begin{csbabstract}
{\noindent{{\bf Abstract}}\\
We introduce the idea that resampling from past observations in a Markov Chain Monte Carlo sampler can fasten convergence. We prove that proper resampling from the past does not disturb the limit distribution of the algorithm. We illustrate the method with two examples. The first on a Bayesian analysis of stochastic volatility models and the other on Bayesian phylogeny reconstruction.}
{Key words: Monte Carlo methods, Resampling, Stochastic volatility models, Bayesian phylogeny reconstruction}
\end{csbabstract}
\section{Introduction}\label{intro}
Markov Chain Monte Carlo (MCMC) methods have become the standard computational tool for bayesian inference. But the great flexibility of the method comes with a price. Namely, it is very difficult to determine a priori (before the simulation) or a posteriori whether a given MCMC sampler can mix or has mixed in a given computing time. The challenge becomes that of designing fast converging Monte Carlo algorithms. Contributions in this field can have significant impact in other scientific disciplines where these methods are used. 

In this paper, we propose a new and general approach to increase the convergence rate of MCMC algorithms. The method is based on resampling. Suppose that at time $n$, we want to sample $X_n$ in a MCMC algorithm. Instead of sampling $X_n$ from $P(X_{n-1},\cdot)$ for some transition kernel $P$, we propose to obtain $X_n$ by resampling independently from $\{X_B,\ldots,X_{n-1}\}$, where $B\geq 0$ is some burn-in period. This \textit{resampling from the past} step is then repeated during the simulation at some predetermined times  $a_1<a_2<\ldots$. Basically, the idea is to look at $\{X_B,\ldots,X_{n-1}\}$ as a sample from $\pi$. Therefore resampling from the past allows the sampler to move more easily and according to a distribution that is close to $\pi$. The resampling schedule plays an important role. As long as we do not resample too much (typically, we need $(a_n)$ such that $a_n/n\to\infty$ as $n\to\infty$), we show that resampling from the past does not disturb the limit distribution of the sampler.

Resampling from the past can perform poorly if the original sampler has a very poor convergence rate. We extend the framework above by allowing resampling from an auxiliary process $\{X^{(0)}_n\}$ that has a better convergence rate towards its target distribution $\pi^{(0)}$. Resampling from an auxiliary process is not new and is the idea behind the equi-energy sampler  recently proposed by \citep{kzw06}. But the equi-energy sampler has a number of complications that we avoid here by using an \textit{importance-resampling}. The idea is also apparent in the ``Metropolis with an adaptive proposal'' of \citep{chauveauetvand01}. On the theoretical side, we show in the case of \textit{importance-resampling}, that resampling from an auxiliary process does not disturb the limit distribution of the sampler.

We apply our methods to two examples from Bayesian data analysis. First, we consider the Bayesian analysis of stochastic volatility models \citep{kimetal97}. We improve the efficiency of the basic Gibbs sampler for this problem by a factor of fifty (50). In the second example, we look at Bayesian phylogenetic trees reconstruction. Our methods improve the efficiency of the MCMC sampler of \citep{largetetsimon99} by a factor of hundred (100).

The paper is organized as follows. In Section \ref{past}, we present the idea of resampling from the past. Resampling from an auxiliary process is discussed in Section \ref{auxil}. All the theoretical proofs are postponed to Section \ref{proof} and the simulation examples are presented in Section \ref{examples}.

\section{Resampling from the past}\label{past}
Let $\{X_n\}$ be a Markov chain with state space $(\X,\B)$, transition kernel $P$ and invariant distribution $\pi$ started at $X_0=x$. If the chain is ergodic then $\mathcal{L}_x(X_n)$, the distribution of $X_n$, will converge to $\pi$ as $n\to\infty$. But it is well known that for MCMC algorithms, the convergence of $\mathcal{L}_x(X_n)$ to $\pi$ can be too slow for the sampler to be useful. We propose the following idea to accelerate the convergence of Markov chains. Suppose that after a burn-in period $B$, we have the sample $\{X_B,X_{B+1},\ldots,X_{n-1}\}$ at time $n$. Instead of sampling $X_n\sim P(X_{n-1},\cdot)$ as we normally do, we obtain $X_n$ by resampling independently and with equal weight from $\{X_B,X_{B+1},\ldots,X_{n-1}\}$. The resampling step is then repeated at some predetermined times $a_1<a_2<\ldots$. Intuitively, if $P$ mixes reasonably well, $\{X_B,X_{B+1},\ldots,X_{n-1}\}$ can be seen as a sample points from $\pi$ and resampling will operate as an i.i.d. sampling from $\pi$. 

Consider the following toy example. We want to use the Random Walk Metropolis (RWM) algorithm with proposal density $q(x,y)=\mathcal{N}(y-x;0,\sigma^2)$ with $\sigma=0.1$ to sample from the standard normal density $\mathcal{N}(x;0,1)$; where $\mathcal{N}(x;\mu,\sigma^2)$ denotes the density of the normal distribution $N(\mu,\sigma^2)$ with mean $\mu$ and variance $\sigma^2$. We compare the plain RWM with a RWM with resampling. Each sampler is run for $25,000$ iterations. Graph 1 (a) shows the last $5,000$ sample points and Graph (b), the autocorrelation function from the last $20,000$ points in the plain RWM sampler. For the RWM with resampling, we resample at times $B+\lceil k^\alpha\rceil$ (see the justification below), with $B=5,000$ and $\alpha=1.3$. Graph 1 (c) and (d) show the corresponding results for the RWM with resampling. As we can see, there is a significant gain in efficiency. 

Intuitively, resampling helps to the extend that $P$ mixes rapidly. Differently put, the slower $P$ converges to $\pi$, the longer we should wait between two resampling. What should be the resampling schedule $(a_k)$? Obviously, we should not resample all the time. We find that the choice $a_k=b_1+b_2k^\alpha$, $\alpha>1$ is a valid choice and works well in practice for $b_2=1$, and $\alpha\approx 1.3$. The choice $a_k=b_1+b_2k$ is also theoretically valid as long as $b_2$, the time between two resampling, is large enough. 
\bigskip

\begin{center}
\scalebox{0.5}{\includegraphics{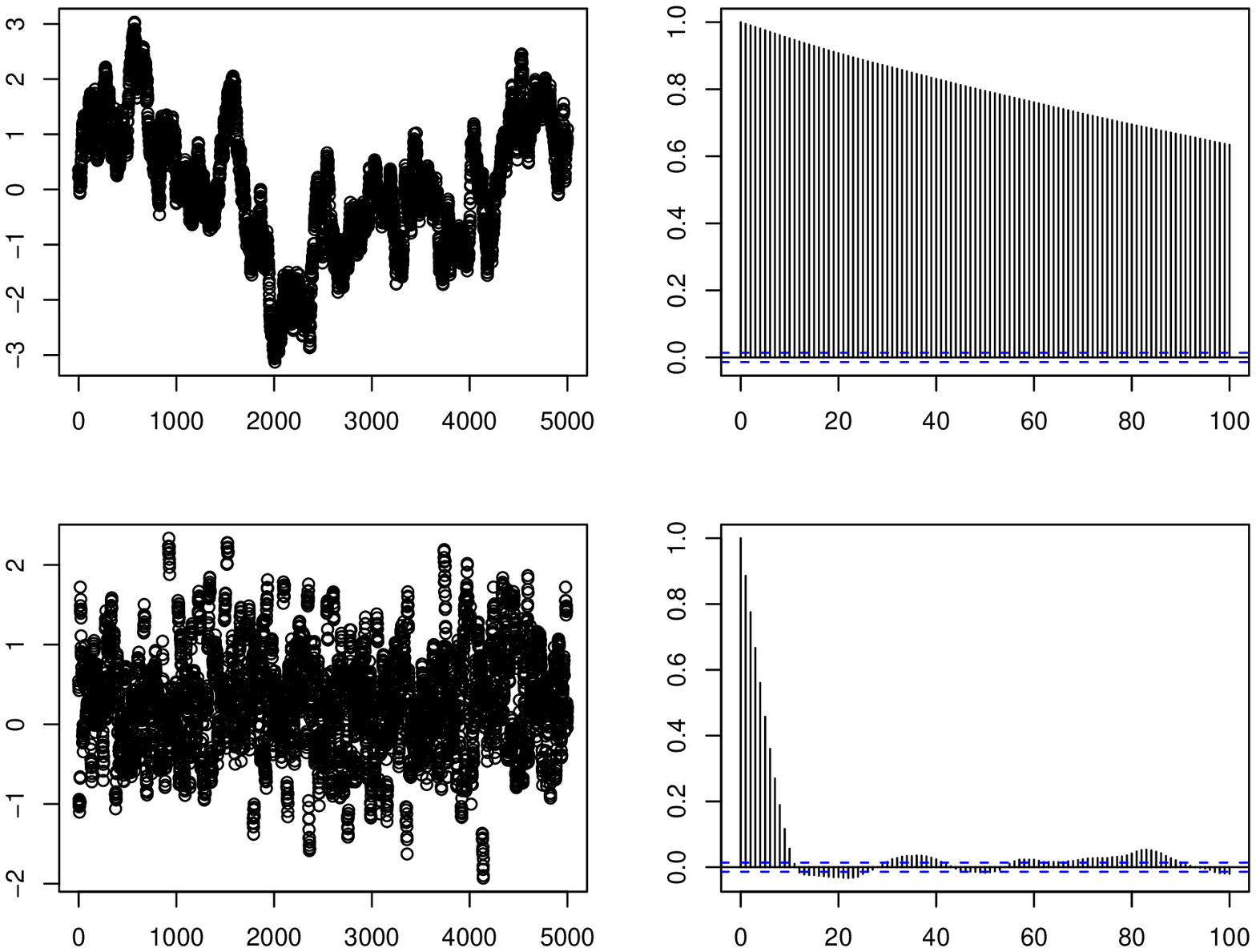}}\\
\noindent\underline{Graph 1}: Comparing a plain RWM and a RWM with resampling in sampling from the standard norma distribution $N(0,1)$.
\end{center}
\bigskip

\subsection{Theoretical discussion}\label{theorie1}
What can we prove about this algorithm? We can prove that despite the resampling, the limit distribution of the algorithm is $\pi$ under certain conditions on $P$ and on the resampling schedule $(a_k)$. We recall the algorithm. The resampling schedule $0<a_1<a_2<\cdots<a_n<\infty$ is given and is nonrandom. Fix $B$ the burn-in period. We start the sampler at some arbitrary point $X_0=x$. At time $n\geq 1$, given $\{X_0,\ldots,X_{n-1}\}$, if $n>B$ and $n=a_k$ for some $k\geq 1$ then $X_n\sim\frac{1}{n-B}\sum_{j=B}^{n-1}\delta_{X_j}(\cdot)$. Otherwise sample $X_n\sim P(X_{n-1},\cdot)$. We denote $\Pr$ the underlying probablity measure and $\E$ its expectation operator. Here are some standard notations that we use below. If $P_1$ and $P_2$ are two transition
kernels on $\X$, the product $P_1P_2$ denotes the transition kernel
$P_1P_2(x,A):=\int P_1(x,dy)P_2(y,A)$. Recursively, we can 
define $P_1^n$ by $P_1^1=P_1$ and $P_1^n=P_1^{n-1}P_1$. A transition kernel $P_1$ defines a linear operator (also denoted $P_1$) on the space of $\Real$-valued functions on $(\X,\B)$ into itself, by $P_1f(x):=\int P_1(x,dy)f(y)$. If $\mu$ is a signed measure on $(\X,\B)$, we denote $\mu(f):=\int\mu(dx)f(x)$ and we will also write $\mu$ to denote the linear functional on the space of $\Real$-valued functions on $(\X,\B)$ thus induced. Finally, we define $\mu P_1(A):=\int\mu(dx)P_1(x,A)$. Let $V:\;\X\to [1,\infty)$ be given. For $f:\;\X\to\Real$, we define its $V$-norm $\abs{f}_V:=\sup_{x\in\X}\frac{\abs{f(x)}}{V(x)}$ and we introduce the space $L_{V}:=\{f:\;\X\to\Real:\;\; \abs{f}_V<\infty\}$. For a signed measure $\mu$ on $(\X,\B)$ we define its $V$-norm $\norm{\mu}_V:=\sup_{f\in L_V,\;\abs{f}_V\leq 1}\abs{\mu(f)}$. Similarly, for a linear operator $T$ from the space of $\Real$-valued functions on $\X$ into itself, we define $\nnorm{T}_V:=\sup_{f\in L_V,\;\abs{f}_V\leq 1}{\abs{Tf}_V}$. If $\nnorm{T}_V<\infty$, then $T$ defines a bounded linear operator from the Banach space $(L_V,\abs{\cdot}_V)$ into itself. 

We assume that the transition kernel $P$ in the algorithm is geometrically ergodic in the sense that:

\noindent\textbf{Assumption (A):} \textit{$P$ is irreducible, aperiodic and there exists $\rho\in (0,1)$, a measurable function $V:\;\X\To
[1,\infty)$ such that 
\begin{equation}\nnorm{P^{n}-\pi}_V=O\left(\rho^n\right),\end{equation} }

This assumption implies that $\pi(V)<\infty$ and that $\sup_nP^nV^\alpha(x)<\infty$ for any $x\in\X$, $\alpha\in[0,1]$. We refer the reader to \citep{meynettweedie93} for more on geometrically ergodic Markov chains. This is a convenient assumption that is known to hold for many MCMC sampler. 

Define $c:=\frac{1}{1-\rho}$ and $\delta_n:=-a_1\log(\rho)+\sum_{k=2}^n\log(a_{k})-\log(c+a_{k-1})$.
\begin{theorem}\label{thm1}
Assume (A). Then there exists a constant $C\in (0,\infty)$ such that for $a_k\leq n<a_{k+1}$:
\begin{equation}\nnorm{\mathcal{L}^{(n)}-\pi}_V\leq C\rho^{n-k}\exp\left[-\delta_k\right],\end{equation}
where the transition kernel $\mathcal{L}^{(n)}$ is defined by $\mathcal{L}^{(n)}(x,A):=\Pr\left[X_n\in A\vert X_0=x\right]$. In particular if $\delta_n\to\infty$ as $n\to\infty$, the algorithm has limit distribution $\pi$.
\end{theorem}
\begin{proof}See Section (\ref{proof}).\end{proof}

\bigskip

Resampling from the past can sensibly reduce the autocorrelation in the output of a MCMC algorithm. But when the sampler has a very slow mixing time, it might be better to resample from an auxiliary process that has a better mixing time.

\section{Resampling from an auxiliary process}\label{auxil}
As above, $\pi(dx)\propto h(x)\lambda(dx)$ is the probability measure of interest on the measure space $(\X,\B,\lambda)$. We introduce another probability measure $\pi^{(0)}(dx)\propto h^{(0)}(x)\lambda(dx)$ on $(\X,\B,\lambda)$. Let $\{X^{(0)}_n\}$ be a Markov chain with invariant distribution $\pi^{(0)}$ and transition kernel $P^{(0)}$. Let $k:\;\X\times\X\to [0,\infty)$ be a measurable function and $T$ a transition kernel on $(\X,\B)$. Define the transition kernel $Q(x,dy)=\frac{\int\pi^{(0)}(dz)k(x,z)T(z,dy)}{\int \pi^{(0)}(dz)k(x,z)}$. Following \citep{tierney98}, let $S\subseteq\X\times\X$ be such that the probability measures $\pi(dx)Q(x,dy)$ and $\pi(dy)Q(y,dx)$ are mutually absolutely continuous on $S$ and mutually singular on $\X\setminus S$. 

We assume that $\{X^{(0)}_n\}$ converges (reasonably quickly) to $\pi^{(0)}$. Let $P$ be a transition kernel with invariant distribution $\pi$ and $\theta\in[0,1]$. The algorithm works as follows. Given $(X_0^{(0)},\ldots,X_n^{(0)},X_0,\ldots,X_n)$:
\begin{itemize}
\item with probability $\theta$, we sample $X_{n+1}$ from $P(X_n,\cdot)$;
\item with probability $1-\theta$, we propose $Y$ from $R_n(X_n,\cdot)$ where $R_n(x,A)=\frac{\sum_{l=0}^nk(x,X_l^{(0)})T(X^{(0)}_l,A)}{\sum_{l=0}^nk(x,X_l^{(0)})}$. In other words, we resample $Y_1$ from $\{X_0^{(0)},\ldots,X_n^{(0)}\}$ with weights $k(X_n,X_l^{(0)})$ and propose $Y\sim T(Y_1,\cdot)$.

Then we either ``accept'' $Y$ and set $X_{n+1}=Y$ with probability $\alpha(X_n,Y)$, or ``reject'' $Y$ and set $X_{n+1}=X_n$ with probability $1-\alpha(X_n,Y)$, where
\begin{equation}\label{acc}
\alpha(x,y)=\left\{\begin{array}{lc}\min\left[1,\frac{\pi(dy)Q(y,dx)}{\pi(dx)Q(x,dy)}\right] & \mbox{ if } (x,y)\in S\\
0 & \mbox{ otherwise}.\end{array}\right.\end{equation}
\end{itemize}

For $n$ large enough, a sample from $R_n(x,\cdot)$ can be seen as a sample from $Q(x,dy)$ which explain the acceptance probability (\ref{acc}). But the algorithm is not feasible as such because the ratio in (\ref{acc}) cannot be computed in general. The natural choice which simplifies $Q$ is to choose a transition kernel $T$ that is invariant under $\pi$ and $k(x,y)=\omega(y)=h(y)/h^{(0)}(y)$. With this choice, we get $\alpha(x,y)\equiv 1$ on $S$. We call this scheme \textit{importance-sampling} resampling. It is not necessary to choose a complicated transition kernel for $T$. Throughout, we choose $T$ to be the identity transition kernel, $T(x,A)=\textbf{1}_A(x)$ in which case $S=\{(x,y):\; 0<h(x)k(x,y)<\infty\}$.

Another choice for which the acceptance ratio $\alpha(x,y)$ simplifies is $T(x,A)=\textbf{1}_A(x)$ and $k(x,y)=\textbf{1}_{\{D(x)\}}(y)$ where $(D_i)$ is a given partition of $\X$ and $D(x)=D_i$ if $x\in D_i$. This corresponds to the set-up of the equi-energy sampler of \citep{kzw06}. With this choice of $k$, the acceptance probability becomes $\alpha(x,y)=\min\left(1,\frac{\omega(y)}{\omega(x)}\right)$ (and $0$ if $\omega(x)=0$ or $\omega(x)=\infty$). The drawback with this choice is that we have to define the partition $(D_i)$ in the first place and an inadequate partition can result in a high rejection rate for the resampling step.


\bigskip
\begin{algo}[MCMC with Importance-Resampling from an auxiliary process]
\bigskip
At some time $n\geq 1$, given $\left(X_0^{(0)},\ldots,X_n^{(0)},X_0,\ldots,X_n\right)$:
\begin{description}
\item [(i)] With probability $\theta$, sample $X_{n+1}$ from $P(X_n,\cdot)$. Otherwise with probability $1-\theta$ sample $X_{n+1}$ from \[\frac{\sum_{i=0}^n\omega(X_i^{(0)})\delta_{X_i^{(0)}}(\cdot)}{\sum_{i=0}^n\omega(X_i^{(0)})}.\]
\item [(ii)] Sample $X_{n+1}^{(0)}$ from $P^{(0)}(X_n^{(0)},\cdot)$.
\end{description}
\end{algo}
\bigskip

\subsection{Theoretical discussion}\label{theorie2}
We look more closely to $\{X_n\}$ when the importance-resampling scheme is used. \citep{atchadeetliu06} have shown that the limit distribution of the equi-energy sampler is indeed $\pi$ under a number of conditions. We can study the process $\{X_n\}$ along the same line. The assumption we impose are less stronger than in \citep{atchadeetliu06}. We continue with the notations in Section \ref{theorie1}. Essentially we will assume that $P^{(0)}$ is geometrically ergodic and that the weight function satisfies $\abs{\omega}_{V^{\alpha}}<\infty$, for some $\alpha\in[0,1/4)$. Typically $\omega$ is bounded.
\bigskip

\noindent\textbf{Assumption (A0):}\textit{ $P^{(0)}$ is irreducible and aperiodic and there exists $\rho_0\in (0,1)$ such that 
\begin{equation}\nnorm{P^{(0)^n}-\pi^{(0)}}_V=O\left(\rho_0^n\right),\end{equation} where $V$ is as in (A).
}

\bigskip

\begin{theorem}\label{thm2}
Assume that $P$ satisfies (A), $P^{(0)}$ satisfies (A0) and $\abs{\omega}_{V^{\alpha}}<\infty$ for some $\alpha\in[0,1/4)$. Then for any measurable function $f:\;\X\to\Real$ such that $\abs{f}_{V^{\alpha}}<\infty$,
\begin{equation}\E\left[f(X_n)\vert X_0=x\right]\To\pi(f), \; \mbox{ as } n\to\infty\end{equation}
and
\begin{equation}\frac{1}{n}\sum_{i=0}^{n-1}\left(f(X_i)-\pi(f)\right)\stackrel{a.s.}{\To} 0,\; \mbox{ as } n\to\infty.\end{equation}
\end{theorem}
\begin{proof}See Section \ref{proof}.\end{proof}

\section{simulation examples}\label{examples}
We illustrate the methods developed above with two examples from bayesian modelling. In the first example, we consider the Bayesian analysis of stochastic volatility models (\citep{kimetal97}) and in the second example, we look at Bayesian phylogenetic trees reconstruction (\citep{largetetsimon99}).
\subsection{Bayesian analysis of stochastic volatility models}
We consider the Bayesian analysis of the basic stochastic volatility model:
\begin{eqnarray}\label{svmodel}y_t&=&e^{h_t/2}\eps_t,\;\;t=0,\ldots,T\\
h_{t+1}&=&\mu+\phi(h_t-\mu)+\sigma u_t,\;\; t=0,\ldots,T-1,\end{eqnarray}
where $(\eps_t)$ and $(u_t)$ are two uncorrelated sequences of i.i.d. standard normal random variables. We assume that $h_0\sim N\left(\mu,\frac{\sigma^2}{1-\phi^2}\right)$ and $\abs{\phi}<1$ to assure the stationarity of the process $(h_t)$. We observe $(y_t)$ but not $(h_t)$, the so-called volatility process. The objective is to estimate $\theta=(\sigma,\phi,\beta)$ where $\beta=e^{\mu/2}$. This model and its generalizations  have attracted attention in the financial econometrics literature as a better way to model financial markets series. A bayesian approach to analyze this model has been proposed by a number of authors (see e.g. \citep{kimetal97} and the references therein). The difficulty is that the volatility process $(h_t)$ is not observed making the likelihood of $\theta$ analytically intractable. The natural solution is to see $(h_t)$ as a parameter and to design a Gibbs sampler on the posterior distribution $\pi(\theta,h_0,\ldots,h_T)$, of the parameter $\theta$ and the volatility process $(h_0,\ldots,h_T)$. But, due to the high autocorrelation in the volatility process, this sampler mixes very slowly. This mixing problem has motivated some authors to propose more sophisticated reparametrization of the model for better MCMC convergence. We show here that by resampling from the past in the Gibbs sampler, we can match the performances of the sophisticated solution proposed in \citep{kimetal97}.

We use the same prior distribution for $\theta$ as in \citep{kimetal97} and essentially the same Gibbs sampler to sample from $\pi(\theta,h_0,\ldots,h_T)$ except when sampling from the conditional $\pi(h_t\vert \theta,h_{-t})$. To sample from this conditional, we use an Independent Metropolis sampler instead of the Accept-Reject method adopted in \citep{kimetal97}. The proposal distribution of our Independent Metropolis sampler is the same as the dominating distribution in the Accept-Reject sampler of \citep{kimetal97}. We refer the reader to \citep{kimetal97} for the details.

Following \citep{kimetal97} and \citep{shephardetpitt97}, we use model (\ref{svmodel}) to analyze the Sterling dataset, which gives the daily observations of weekday close exchange rates for the UK Sterling/US Dollar exchange rate from $1/10/81$ to $28/6/85$. The total number of observations is $T=946$. We first center the series with the formula \\$y_t=100\left[(\log(r_t)-\log(r_{t-1}))-\frac{1}{n}\sum_{j=1}^n(\log(r_j)-\log(r_{j-1}))\right]$, where $(r_t)$ is the observed exchange rates. We then model $(y_t)$ with the model (\ref{svmodel}).

We compare the plain Gibbs sampler with the 2 strategies discussed above: a Gibbs sampler with resampling from the past and a Gibbs sampler with resampling from an auxiliary process. To assure that the three sampler have about the same computational cost (storage requirement aside), we set the auxiliary process to be another copy of the plain Gibbs sampler with the same target distribution. The three samplers are run for $N=250,000$ iterations. For each sampler and for each of the variables $\sigma$, $\phi$, $\beta$, we give a plot of the last $5,000$ sample points together with the histogram and the autocorrelation function from the last $100,000$ points. When resampling from the past, the resampling schedule used is $B+\lceil k\rceil^\alpha$, $B=125,000$ and $\alpha=1.25$. For the third sampler with resampling from an auxiliary process, each of the two chains is run for $125,000$ iterations. The results of the variable $\sigma$ (resp. $\phi$ and $\beta$) are given in in Graph 2 (resp. Graph 3 and Graph 4). On each graphics, the first column gives the result of the plain Gibbs sampler, the second column gives the results of the Gibbs sampler with resampling from the past and the results of the third sampler are in the third column. 

Clearly, resampling from the past significantly improve on the Gibbs sampler. To quantify the gain, we compute, following \citep{kimetal97} the inefficiency of each sampler on each of the three variables. For a Markov chain with transition kernel $P$ and invariant distribution $\pi$, the inefficiency at $f$ is:
\begin{equation}
I(f)=1+2\sum_{k=1}^\infty\rho_k(f),\end{equation}
where $\rho_k(f)=Cov_\pi\left(f(X_k),f(X_0)\right)/Var_\pi\left(f(X_0)\right)=\pi\left(\bar fP^k\bar f\right)/\pi\left(\bar f^2\right)$.  Basically, it is the cost of using a dependent process to sample from $\pi$. To estimate $I(f)$, we use, following \cite{kimetal97}: 
\begin{equation}
\hat I(f)=1+\frac{2B}{B-1}\sum_{i=1}^BK\left(\frac{i}{B}\right)\hat\rho_i(f),\end{equation}
where $\hat\rho_i(f)$ is the usual estimate of the autocorrelation at lag $i$ for $f$ and $K$ the so-called Parzen kernel. We use $B=5,000$. The result is given in Table 1.
 
By resampling from the past or from an auxiliary process, we obtain a sampler that outperforms \citep{shephardetpitt97} and is as efficient as the offset mixture method of \citep{kimetal97}.

\begin{table}
\begin{center}
\small
\begin{tabular}{lccc}
\hline
&$\sigma$&$\phi$&$\beta$\\
\hline
Plain Gibbs& $448.12$&$211.55$ & $1.54$\\
Gibbs with resampling & $10.96$& $4.92$ & $0.97$\\
Gibbs with Aux. Proc. & $12.24$& $9.91$ & $1.39$\\
\hline
\end{tabular}
\end{center}
\caption{Inefficiencies of the samplers for the Sterling dataset.}
\end{table}

\subsection{Bayesian phylogeny reconstruction}
Since Darwin's theory of evolution, methods to reconstruct the evolutionary relationships between different species have become important. We are concerned here with the statistical inference of phylogenetic trees based on molecular sequences. Recently, more realistic models have been considered in this field owing to the MCMC machinery. We show here that MCMC samplers for phylogeny reconstruction can be improved upon with resampling from the past.

The statistical model is not standard, so we summarize it first. For more details on phylogenetic trees, we refer the reader to \citep{felsenstein04}. Suppose we have $n$ aligned deoxyribonucleic acid (DNA) sequences $(y_1,\ldots,y_n)$ each of length $m$, where sequence $i$ is from organism $i$. That is, $y_i=(y_i(1),\ldots,y_i(m))$ where $y_i(j)$ can be one of the four nucleotide basis $A$ (Adenine), $G$ (Guanine), $C$ (Cytosine) or $T$ (Thymine). Based on these sequences, we would like to infere the phylogenetic tree or evolutionary relationships between these organisms. To be precise, we recall that a binary tree $\tau$ for $n$ species is a connected graph $(V,E)$ with vertex set $V$ and edges $E$, with no cycle, such that $V=\{\rho\}\cup\mathcal{I}\cup\mathcal{T}$, where $\rho$ (the root) has degre $2$; any $v\in\mathcal{I}$ has degre $3$ and any $v\in\mathcal{T}$ has degre $1$. $\mathcal{I}$ has $n-2$ elements called the internal nodes and $\mathcal{T}$ (the leaves or the tips) represent the $n$ species. A phylogenetic tree for $n$ species is a couple $\psi=(\tau,b)$, where $\tau$ is a binary tree for the $n$ species and $b\in (0,\infty)^{\abs{E}}$, where $\abs{E}=2n-1$ is the cardinality of $E$. For $e\in E$, $b_e$ represents the length of edge $e$, the so-called branch length. We restrict our attention to phylogenetic trees with ``contemporary tips'', where the sum of the branch length $b_e$ on the directed path from the root to any tip is constant (equal to $1$ hereafter). Such phylogenetic trees are said to be with a ``molecular clock'' as the $b_e$ can now be interpreted as time. Let $\Psi$ be the set of all phylogenetic trees for $n$ species. For $i\in V\setminus\{\rho\}$, denote $p(i)$ the parent of $i$, that is the vertex $p(i)$ such that $(p(i),i)\in E$.

The model of phylogenetic reconstruction we are interested in assumes that there are some missing DNA sequences $(y_j)_{\{j\in\{\rho\}\cup\mathcal{I}\}}$ such that the joint conditional distribution of $(y_j)_{V}$ given the phylogenetic tree $\psi$ writes:
\begin{equation}\label{phylomodel1}
f((y_i)_{\{i\in V\}}\vert\psi)=f(y_\rho)\prod_{i\in V\setminus\{\rho\}}f\left(y_i\vert y_{p(i)},\psi\right).\end{equation}

In addition we make the simplifying assumption that each site evolves independently:
\begin{eqnarray}\label{phylomodel2}
f(y_\rho)&=&\prod_{j=1}^mf(y_\rho(j)),\;\mbox{ and }\\
f\left(y_i\vert y_{p(i)},\psi\right)&=&\prod_{j=1}^m f\left(y_i(j)\vert y_{p(i)}(j),b_{(p(i),i)}\right).\end{eqnarray}

And finally, we assume that there exist $(\pi_l)_{l\in \{A,G,C,T\}}$, $\pi_l\geq 0$, $\sum\pi_l=1$, parameters $\theta,\kappa\in(0,\infty)$ and a $4\times 4$ Markov process generator $Q=Q(\theta,\kappa,\pi_A,\pi_G,\pi_C,\pi_T)$ such that:
\begin{eqnarray}\label{phylomodel3}
f(y_\rho(j)=l)&=&\pi_l,\;l\in\{A,G,C,T\}\mbox{ and }\\
f\left(y_i(j)=m\vert y_{p(i)}(j)=l,b_{(p(i),i)}=b\right)&=&exp(bQ)_{lm},\;\;l,m\in\{A,G,C,T\}.\end{eqnarray}

The matrix $Q$ specifies the model of DNA evolution. We use the F84 model as in \citep{largetetsimon99}. The parameters of the statistical model are then $(\psi,\theta,\kappa,\pi_A,\pi_G,\pi_C,\pi_T)$. To simplify the sampler, we fix $\pi_A,\pi_G,\pi_C,\pi_T$ to their empirical values in the data. We assume that $\psi$ has a uniform prior distribution on $\Psi$ and we assume that $\theta$ and $\kappa$ each has a uniform prior on $(0,M)$, $M=200$. Let $\pi\left(\psi,\theta,\kappa\vert (y)_{i\in\mathcal{T}}\right)$ be the posterior distribution of the model. Clearly, $\pi\left(\psi,\theta,\kappa\vert (y)_{i\in\mathcal{T}}\right)\propto f\left((y)_{i\in\mathcal{T}}\vert \psi,\theta,\kappa\right)$ and this likelihood is obtained by integrating out the missing variables $(y_i)_{i\in\{\rho\}\cup\mathcal{I}}$ from (\ref{phylomodel1}). A fast computation of this likelihood is available with the pruning method of Felsenstein \citep{felsenstein04}. To sample from this posterior distribution, we follow essentially \citep{largetetsimon99}. We update $\theta$ and $\kappa$ together, given the phylogenetic tree $\phi$, using a random walk Metropolis move. Next, given $\theta,\kappa$, we update the phylogenetic tree $\psi$ with the \textit{global move with a molecular clock} of \citep{largetetsimon99}.

We compare this plain MCMC sampler with the samplers obtained with the two methods discussed in this paper. For the simulations, we use the \textit{primate dataset} discussed in \citep{yangetrannala97}. The dataset has $n=9$ species and the phylogeny reconstruction is based on aligned sequences of length $m=888$. The three samplers are simulated for $N=500,000$ iterations. For each sampler and for each of the variables $\theta$, $\kappa$, we give a plot of the last $5,000$ sample points together with the histogram and the autocorrelation function from the last $150,000$ iterations. When resampling from the past, the resampling schedule used is $B+\lceil k\rceil^\alpha$, $B=100,000$ and $\alpha=1.3$. For the third sampler with resampling from an auxiliary process, each of the two chains is run for $250,000$ iterations. The auxiliary process is a MCMC chain with stationary distribution $\pi^{(0)}=\pi^{1/T}$, with $T=2$. The results of the variable $\theta$ (resp. $\kappa$) are given in in Graph 5 (resp. Graph 6). On each graphics, the first column gives the result of the plain MCMC sampler, the second column gives the results of the MCMC sampler with resampling from the past and the results of the third sampler are in the third column. In accordance with \citep{largetetsimon99}, the outputs of the three samplers overwhelmingly (with an estimated posterior distribution over $0.95$) select the phylogenetic tree topology plotted in figure 7 as the most probable for this primate dataset. 

Here again, resampling from the past significantly improve on the plain MCMC sampler. Table 2 gives the efficiency gains.

\begin{table}
\begin{center}
\small
\begin{tabular}{lcc}
\hline
&$\theta$&$\kappa$\\
\hline
Plain MCMC& $1510.23$&$1271.87$\\
MCMC with resampling & $13.95$& $24.37$\\
MCMC with Aux. Proc. & $9.18$& $8.15$\\
\hline
\end{tabular}
\end{center}
\caption{Inefficiencies of the samplers for the primates dataset}
\end{table}
\bigskip

\section{Proofs of Theorem \ref{thm1} and \ref{thm2}}\label{proof}
We start with Theorem \ref{thm1}. Without any loss of generality we assume that $B$, the burn-in period is $0$.
\subsection{Proof of Theorem \ref{thm1}}
The following lemma is a consequence of (A).
\begin{lemma}\label{lem1}
Assume (A). There exists a constant $C_1\in (0,\infty)$ such that for any signed measure $\mu$ on $(\X,\B)$ such that $\mu(\X)=0$ and for any $n\geq 0$,
\begin{equation}\norm{\mu P^n}_V\leq C_1\rho^n\norm{\mu}_V.\end{equation}
\end{lemma}

\begin{proof}[Proof of Theorem \ref{thm1}] Fix $n$ such that $a_k\leq n<a_{k+1}$, $k\geq 2$. For $f\in L_V$ such that $\abs{f}_V\leq 1$, define $\bar f=f-\pi(f)$. We have:
\begin{eqnarray}\E\left(\bar f(X_n)\vert X_0=x\right)&=&\E\left[\E\left(f(X_n)\vert X_{a_k}\right)\vert X_0=x\right]\nonumber\\
&=&\E\left(P^{n-a_k}\bar f(X_{a_k})\vert X_0=x\right)\\
&=&\left(\mathcal{L}^{(a_k)}-\pi\right)\left[P^{n-a_k}\bar f\,\right](x),\end{eqnarray}
where $\mathcal{L}^{(a_k)}(x,A)=\Pr\left(X_{a_k}\in A\vert X_0=x\right)$.  Therefore, since\\ $\nnorm{\mathcal{L}^{(n)}-\pi}_V=\sup_{x\in\X}\frac{\sup_{\abs{f}_V\leq 1}\abs{\E\left(\bar f(X_n)\vert X_0=x\right)}}{V(x)}$, it follows from Lemma \ref{lem1}, that:
\begin{equation}\label{bound1proofthm1}\nnorm{\mathcal{L}^{(n)}-\pi}_V\leq C_1\rho^{n-a_k}\nnorm{\mathcal{L}^{(a_k)}-\pi}_V.\end{equation}

Also, for $f\in L_V$ with $\abs{f}_V\leq 1$, we have:
\begin{eqnarray}
\mathcal{L}^{(a_k)}\bar f(x)&=&\E\left[\frac{1}{a_k}\sum_{j=0}^{a_k-1}\bar f(X_j)\vert X_0=x\right]\nonumber\\
&=&\frac{a_{k-1}}{a_k}\mathcal{L}^{(a_{k-1})}\bar f(x)+\frac{1}{a_k}\sum_{j=a_{k-1}}^{a_k-1}\E\left(\bar f(X_j)\vert X_0=x\right)\nonumber\\
&=&\frac{a_{k-1}}{a_k}\mathcal{L}^{(a_{k-1})}\bar f(x)+\frac{1}{a_k}\sum_{j=a_{k-1}}^{a_k-1}\left(\mathcal{L}^{(a_{k-1})}-\pi\right)P^{j-a_{k-1}}\bar f(x).\end{eqnarray}

Then proceding as above and using Lemma \ref{lem1} again we get:
\begin{equation}\label{boundproof2thm1}\nnorm{\mathcal{L}^{(a_k)}-\pi}_V\leq \exp\left(-u_{k}\right)\nnorm{\mathcal{L}^{(a_{k-1})}-\pi}_V,\end{equation}
with $u_k=\log(a_k)-\log(a_{k-1}+c)$, $c=\frac{1}{1-\rho}$. If we define $u_1=-a_1\log(\rho)$ and $\delta_k=\sum_{i=1}^ku_k$, we get $\nnorm{\mathcal{L}^{(a_k)}-\pi}_V\leq C_2\exp(-\delta_k)$ for some finite constant $C_2$, which, together with (\ref{bound1proofthm1}) yields: 
\begin{equation}\nnorm{\mathcal{L}^{(n)}-\pi}_V=O\left(\rho^{n-a_k}\exp(-\delta_k)\right),\end{equation}
for $a_k\leq n<a_{k+1}$, as wanted.

\end{proof}

\subsection{Proof of Theorem \ref{thm2}}
Let $\{X_n\}$ be the process generated by the importance-resampling scheme.  We prove Theorem \ref{thm2} as a consequence of Theorems 3.1 and 3.2 of \citep{atchadeetrosenthal03}. Denote $\F_n$ the $\sigma$-algebra generated by $(X_0,\ldots,X_n)$. For $x\in\X$ and $A\in\B$, define $P_n(x,A)=\Pr\left(X_n\in A\vert X_{n-1}=x\right)=\Pr\left(X_n\in A\vert \F_{n-1},X_{n-1}=x\right)$. We have:
\begin{equation}
P_n(x,A)=\theta P(x,A)+(1-\theta)\mu_n(A),\end{equation}
where $\mu_n(A)=\E\left[\frac{\sum_{k=0}^{n-1}\omega(X_k^{(0)})\textbf{1}_A(X_k^{(0)})}{\sum_{j=0}^{n-1}\omega(X_j^{(0)})}\right]$.

Define $M_r=\sup_n\E\left(V^{r}(X_n^{(0)})\right)$, $r\geq 0$. It follows from (A0) that $M_r\leq M_1<\infty$ for all $r\in [0,1]$. For $p\geq 0$, we write $\omega_i=\omega(X_i^{(0)})$, $s_n=\sum_{k=0}^{n-1}\omega_k$, $V_i^\alpha=V^{\alpha}(X^{(0)}_i)$ and $\mu_n^{(p)}=\E\left[\frac{\sum_{i=0}^{n-1}\omega_iV_i^\alpha}{s_n}\right]^p$. The next lemma is crutial.
\begin{lemma}\label{lem21} For $p\in[1,\frac{1}{4\alpha}]$, $\max_{0\leq i\leq n-1}\E\left[\frac{\omega_iV_i^\alpha}{s_n}\right]^p=O\left(\frac{1}{n^p}\right)$ and $\mu_n^{(p)}=O(1)$ as $n\to\infty$.\end{lemma}
\begin{proof}
By the Minkowski inequality, we only need to prove that $\max_{0\leq i\leq n-1}\E\left[\frac{\omega_iV_i^\alpha}{s_n}\right]^p=O\left(\frac{1}{n^p}\right)$.

Write $c=\lambda(h)$ and $c_0=\lambda(h^{(0)})$. For $0\leq i\leq n-1$ and $\kappa\in(0,c/c_0)$, we have:
\begin{equation*}\E\left[\frac{\omega_iV_i^\alpha}{s_n}\right]^p=\E\left[\frac{\omega_iV_i^\alpha}{s_n}\textbf{1}_{\{s_n\geq n(c/c_0-\kappa)\}}\right]^p+\E\left[\frac{\omega_iV_i^\alpha}{s_n}\textbf{1}_{\{s_n<n(c/c_0-\kappa)\}}\right]^p.\end{equation*}

\[\E\left[\frac{\omega_iV_i^\alpha}{s_n}\textbf{1}_{\{s_n\geq n(c/c_0-\kappa)\}}\right]^p\leq \frac{1}{n^p(c/c_0-\kappa)^p}\E\left[\omega_iV_i^\alpha\right]^p\leq \frac{\abs{\omega}_{V^{\alpha}}M_1}{n^p(c/c_0-\kappa)^p}.\]

By the Cauchy-Schwarz inequality, we can bound the second term as follows:
\begin{eqnarray*}
\E\left[\frac{\omega_iV_i^\alpha}{s_n}\textbf{1}_{\{s_n<n(c/c_0-\kappa)\}}\right]^p&\leq&\E^{1/2}\left[\frac{\omega_iV_i^\alpha}{s_n}\right]^{2p}\left(\Pr\left[\frac{1}{n}\sum_{i=0}^{n-1}\left(\omega_i-\frac{c}{c_0}\right)<-\kappa\right]\right)^{1/2}\\
&\leq&\E^{1/2}\left[V_i^{2p\alpha}\right]\left(\Pr\left[\frac{1}{n}\sum_{i=0}^{n-1}\left(\omega_i-\frac{c}{c_0}\right)<-\kappa\right]\right)^{1/2}\\
&\leq&\frac{M_1^{1/2}}{n^{2p}\kappa^{4p}}\E^{1/2}\left(\sum_{i=0}^{n-1}\left(\omega_i-\frac{c}{c_0}\right)\right)^{4p},\end{eqnarray*}

where for the last line, the Markov inequality was used. Now we use the classical Poisson equation and martingale approximation technique. Since $\omega\leq V^{\alpha}$, the Poisson equation $\omega-c/c_0=g-P^{(0)}g$ has a solution $g$ which satisfies $\abs{g}\leq V^{\alpha}$. With this solution, for $n>1$, we can rewrite $\sum_{i=0}^{n-1}\omega_i-c/c_0=M_n+W_n$ where $W_n=g(X_0^{(0)})-P^{(0)}g(X_{n-1}^{(0)})$, $M_n=\sum_{i=1}^{n-1}g(X^{(0)}_i)-P^{(0)}g(X^{(0)}_{i-1})$ and $(M_n)$ is a martingale. Therefore with the Minkowski inequality, we get: $\E^{1/2}\left(\sum_{i=0}^{n-1}(\omega_i-c/c_0)\right)^{4p}\leq \left[\E^{1/4p}\left(M_n\right)^{4p}+\E^{1/4p}\left(W_n\right)^{4p}\right]^{2p}$. Since $\abs{g}_{V^\alpha}<\infty$ and $4p\alpha\leq 1$, it follows from Assumption (A0) that $\sup_{i,j}\E\left(g(X_i^{(0)})-P^{(0)}g(X_{j}^{(0)})\right)^{4p}<\infty$. Therefore $\left(\E^{1/4p}\left(W_n\right)^{4p}\right)$ is bounded. Using Burkholder's inequality (see e.g. \citep{halletheyde80}), we have the bound:
\begin{eqnarray*}\E\left(M_n\right)^{4p}&\leq& K_3\E\left(\sum_{i=1}^{n-1}\left(g(X_i^{(0)})-P^{(0)}g(X_{i-1}^{(0)})\right)^2\right)^{2p}\\
&\leq&K_3\left[\sum_{i=1}^{n-1}\E^{1/2p}\left(g(X_i^{(0)})-P^{(0)}g(X_{i-1}^{(0)})\right)^{4p}\right]^{2p}\\
&\leq&K_4n^{2p},\end{eqnarray*}
for some finite constants $K_3,K_4$. This implies that $\E^{1/2}\left(\sum_{i=0}^{n-1}(\omega_i-c/c_0)\right)^{4p}=O(n^p)$ which finishes the proof.

\end{proof}

\begin{lemma}\label{lem2}
For all $n\geq 1$, $P_n$ has an invariant distribution $\pi_n$, and for all $k\geq 0$,
\begin{equation}
\nnorm{P_n^k-\pi_n}_{V^{\alpha}}\leq C\theta^k\rho^k,\end{equation}
where the constant $C\in(0,\infty)$ does not depend on $n$ or $k$. Moreover 
\begin{equation}\pi_n(f)\To\pi(f),\;\mbox{ as } n\to\infty,\end{equation}
for any measurable function $f$, with $\abs{f}_{V^\alpha}<\infty$. \end{lemma} 
\begin{proof}
One can directly check that the invariant distribution of $P_n$ is $\pi_n$ where:
\begin{equation}\pi_n(A)=(1-\theta)\mu_n\left(\sum_{i=0}^\infty \theta^iP^i(x,A)\right).\end{equation}

And by recurrence, we can check that for $k\geq 0$ and $g\in L_{V^{\alpha}}$:
\begin{equation}
P_n^kg-\pi_n(g)=\theta^kP^k\bar g-(1-\theta)\mu_n\left(\sum_{i=k}^\infty\theta^iP^i\bar g\right).\end{equation}

Therefore $\nnorm{P_n^k-\pi_n}_{V^{\alpha}}\leq \theta^k\rho^k\left(1+\frac{1-\theta}{1-\theta\rho}\sup_n\mu_n(V^{\alpha})\right)$ and according to Lemma \ref{lem21},\\ $\sup_n\mu_n(V^{\alpha})$ is finite.

For $f\in L_{V^{\alpha}}$, we write $\zeta(f)=(1-\theta)\sum_{i=0}^\infty\theta^iP^if\in L_{V^{\alpha}}$.  We have $\abs{\pi_n(f)-\pi(f)}=\abs{\mu_n\left(\zeta(\bar f)\right)}$, where $\bar f=f-\pi(f)$. Note that $\pi(\zeta(\bar f))=0$. We recall: 
\begin{equation}\label{eqlem2}\mu_n\left(\bar f\right)=\E\left[\frac{\sum_{k=0}^{n-1}\omega(X_k^{(0)})\bar f(X_k^{(0)})}{\sum_{j=0}^{n-1}\omega(X_j^{(0)})}\right].\end{equation}
 From the strong law of large numbers for $\{X^{(0)}\}$, the expression under the expectation in (\ref{eqlem2}) converges a.s. to $0$ as $n\to\infty$. On the other hand, for $p\in (1,1/4\alpha)$, 
\begin{equation}\E\abs{\frac{\sum_{k=0}^{n-1}\omega(X_k^{(0)})\bar f(X_k^{(0)})}{\sum_{j=0}^{n-1}\omega(X_j^{(0)})}}^p\leq \mu_n^{(p)},\end{equation}
and $(\mu_n^{(p)})$ is a bounded sequence. Therefore the sequence $\left(\frac{\sum_{k=0}^{n-1}\omega(X_k^{(0)})\bar f(X_k^{(0)})}{\sum_{j=0}^{n-1}\omega(X_j^{(0)})}\right)$ is uniformly integrable and it follows that $\mu_n(\bar f)\to 0$ as $n\to\infty$.
\end{proof}

\begin{lemma}\label{lem3}
\begin{equation}\nnorm{P_n-P_{n-1}}_{V^{\alpha}}+\norm{\pi_n-\pi_{n-1}}_{V^{\alpha}}=O\left(\frac{1}{n}\right).\end{equation}
\end{lemma}
\begin{proof}
For $n\geq 1$, we have: $\nnorm{P_n-P_{n-1}}_{V^{\alpha}}+\norm{\pi_n-\pi_{n-1}}_{V^{\alpha}}\leq 2(1-\theta)\E\left[\frac{\omega_{n-1}V_{n-1}^{\alpha}}{\sum_{k=0}^{n-1}\omega_k}\right]$ and the lemma follows from Lemma \ref{lem21}.
\end{proof}
\bigskip

\begin{proof}[Proof of Theorem \ref{thm2}]
Follows from Lemmas \ref{lem2} and \ref{lem3} and Theorems 3.1, 3.2 of \cite{atchadeetrosenthal03}. 
\end{proof}

\bibliographystyle{ims}
\bibliography{biblio}
\newpage

\begin{center}
\scalebox{0.5}{\includegraphics{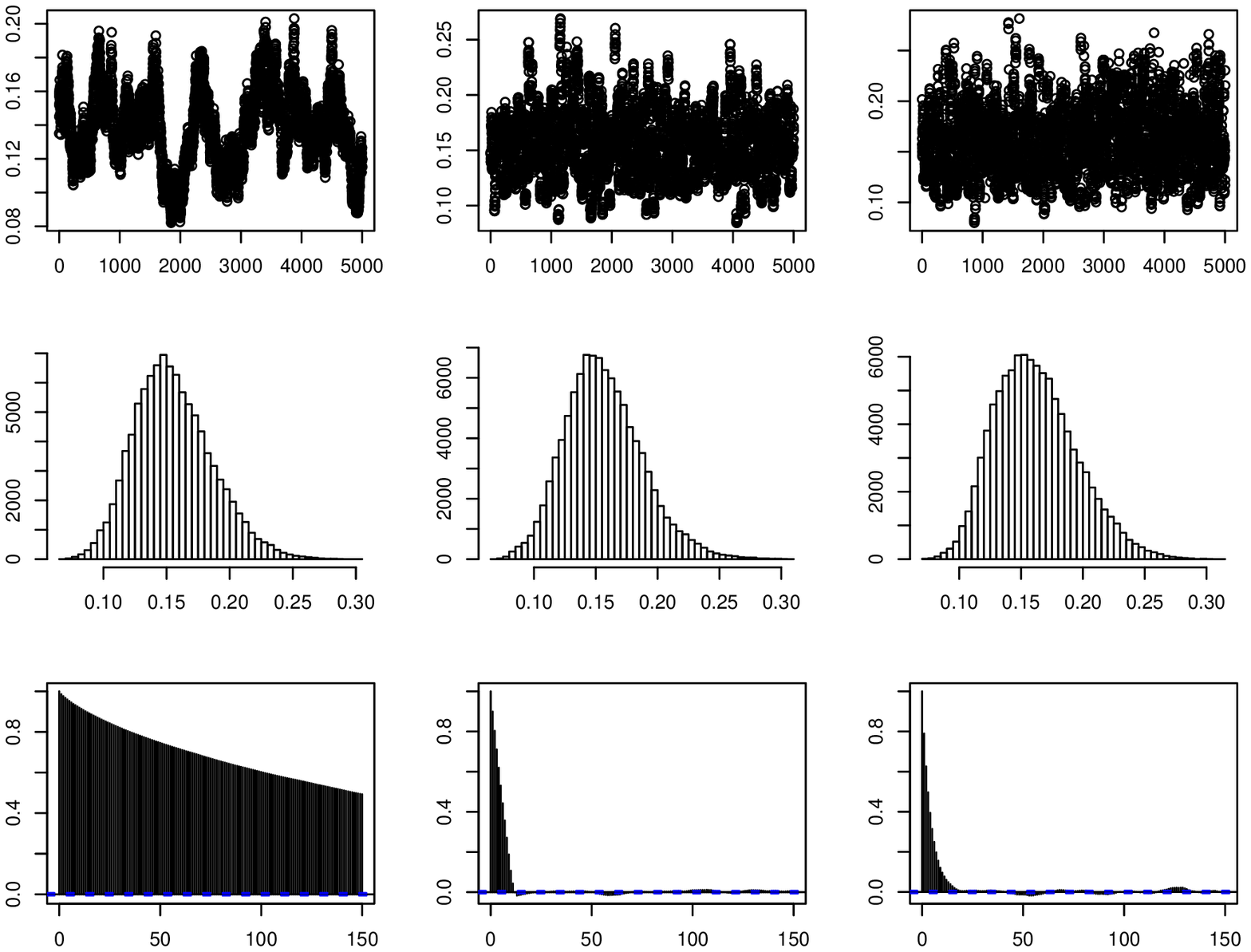}}\\
\noindent\underline{Graph 2}: Outputs for $\sigma$. Sterling dataset. First column is the plain Gibbs, second column is resampling from the past; last column: resampling from an auxiliry Gibbs sampler.
\end{center}

\begin{center}
\scalebox{0.5}{\includegraphics{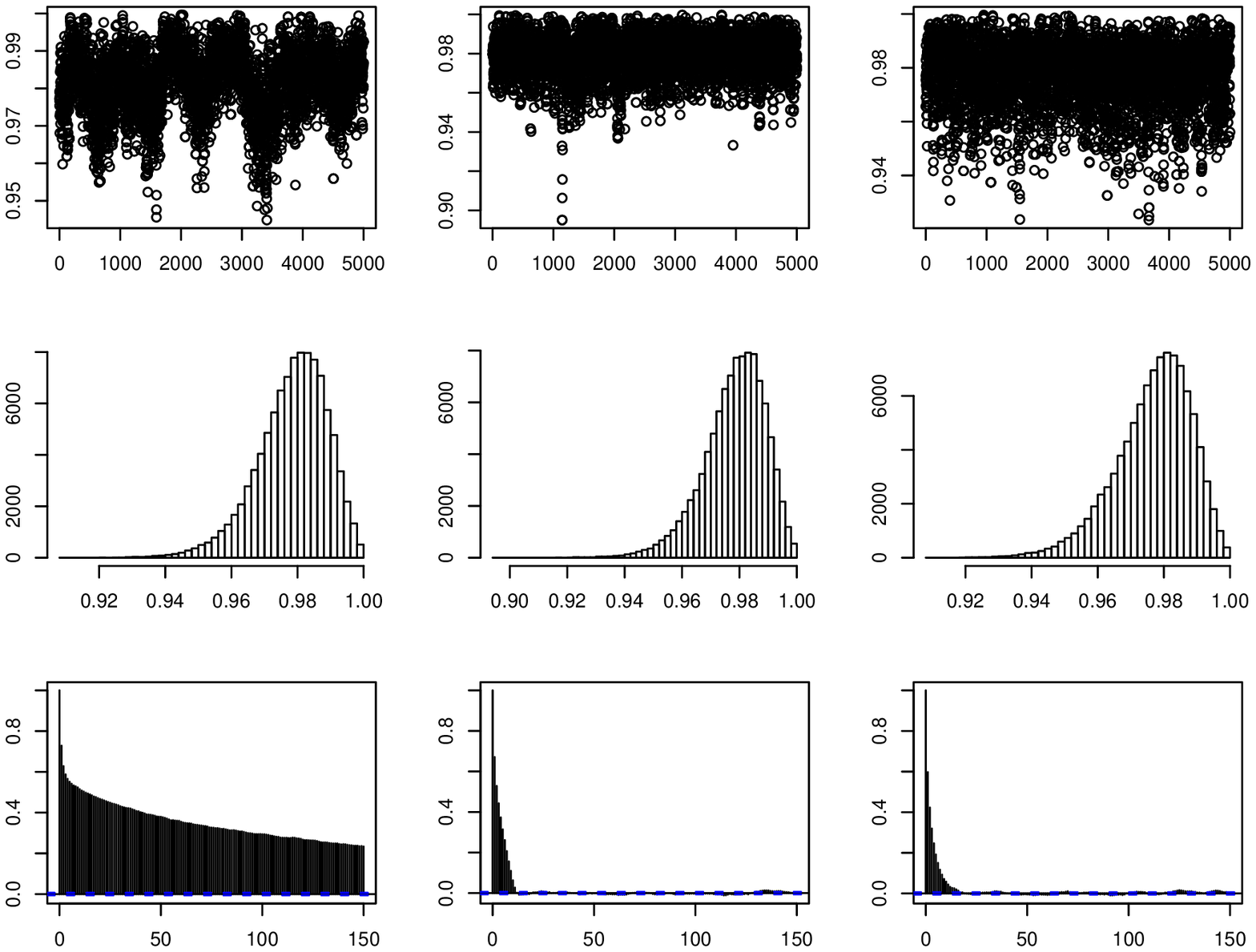}}\\
\noindent\underline{Graph 3}: Outputs for $\phi$. Sterling dataset. First column is the plain Gibbs, second column is resampling from the past; last column: resampling from an auxiliry Gibbs sampler.
\end{center}
\bigskip

\begin{center}
\scalebox{0.5}{\includegraphics{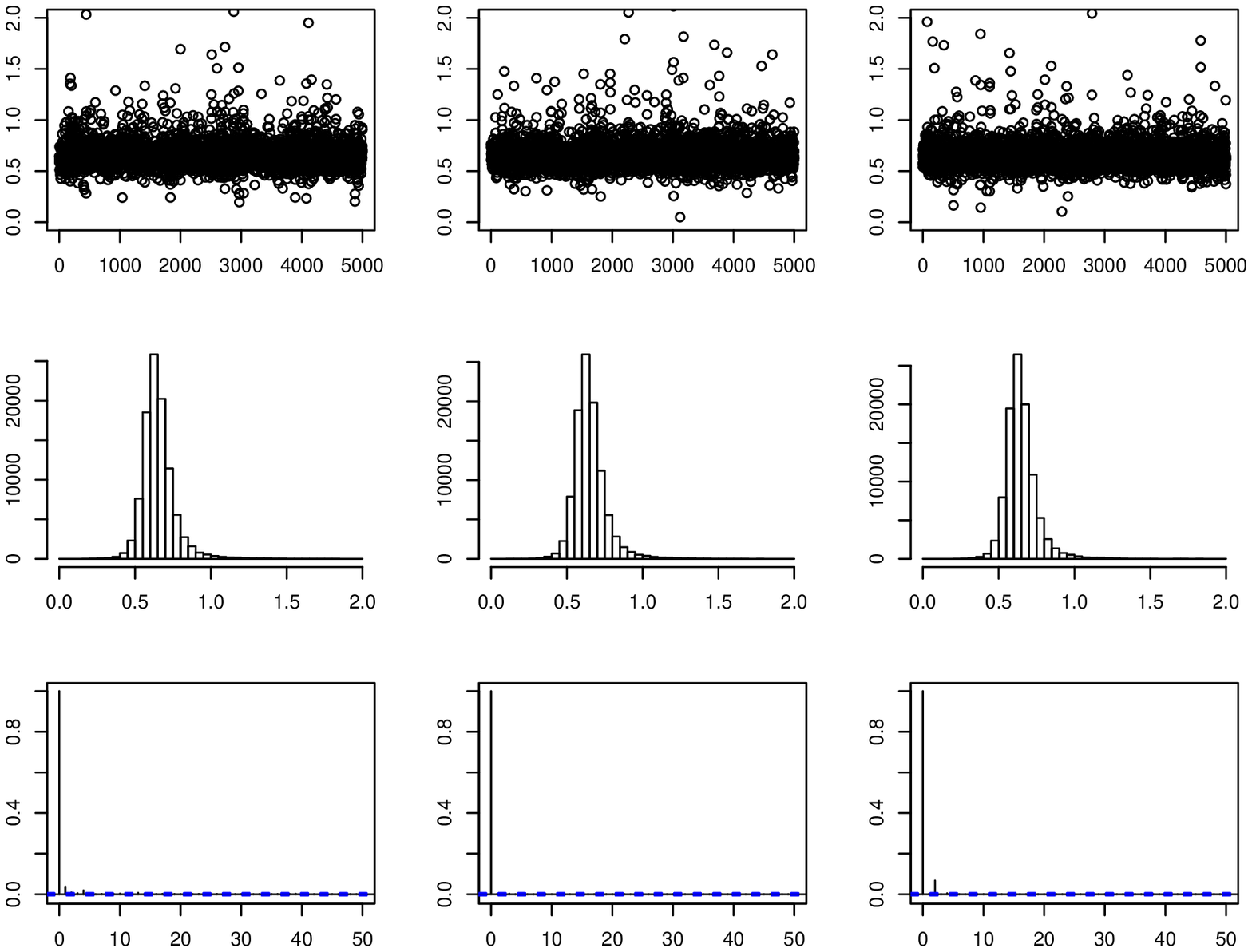}}\\
\noindent\underline{Graph 4}: Outputs for $\beta$. Sterling dataset. First column is the plain Gibbs, second column is resampling from the past; last column: resampling from an auxiliary Gibbs sampler.
\end{center}

\begin{center}
\scalebox{0.5}{\includegraphics{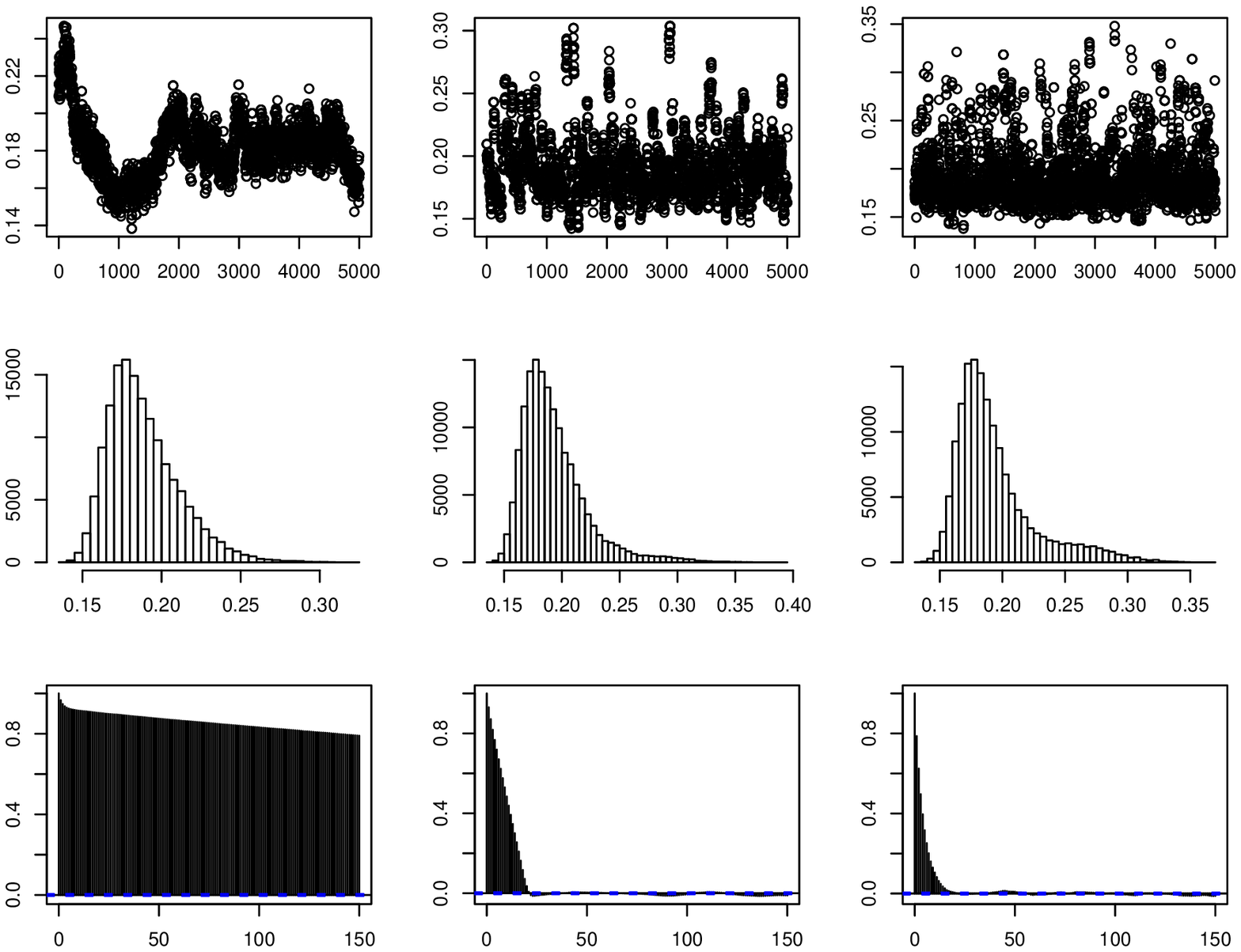}}\\
\noindent\underline{Graph 5}: Outputs for $\theta$. Primates dataset. First column is the plain MCMC, second column is resampling from the past; last column: resampling from an auxiliary MCMC sampler.
\end{center}

\begin{center}
\scalebox{0.5}{\includegraphics{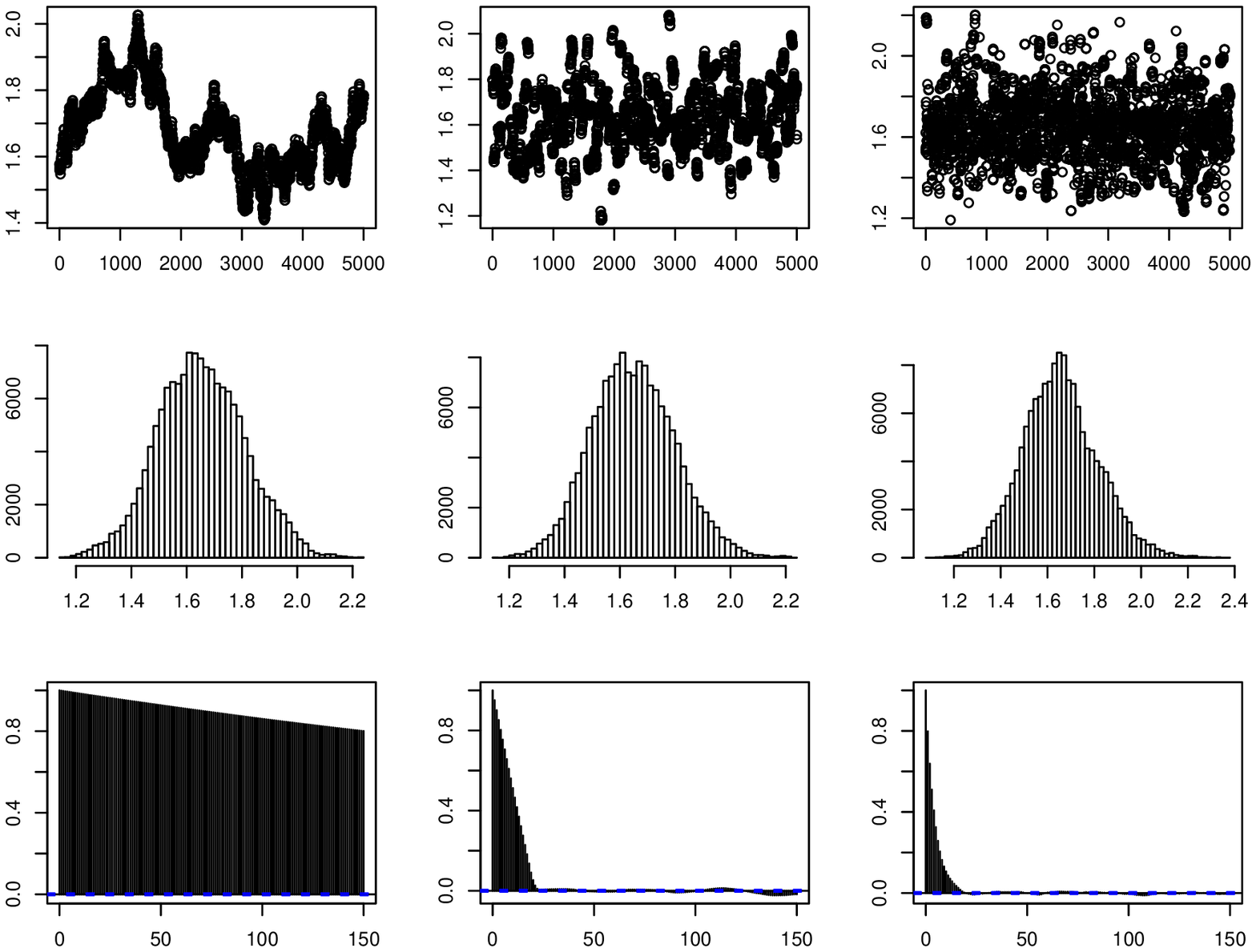}}\\
\noindent\underline{Graph 6}: Outputs for $\kappa$. Primates dataset. First column is the plain MCMC, second column is resampling from the past; last column: resampling from an auxiliary MCMC sampler.
\end{center}

\begin{center}
\scalebox{0.5}{\includegraphics{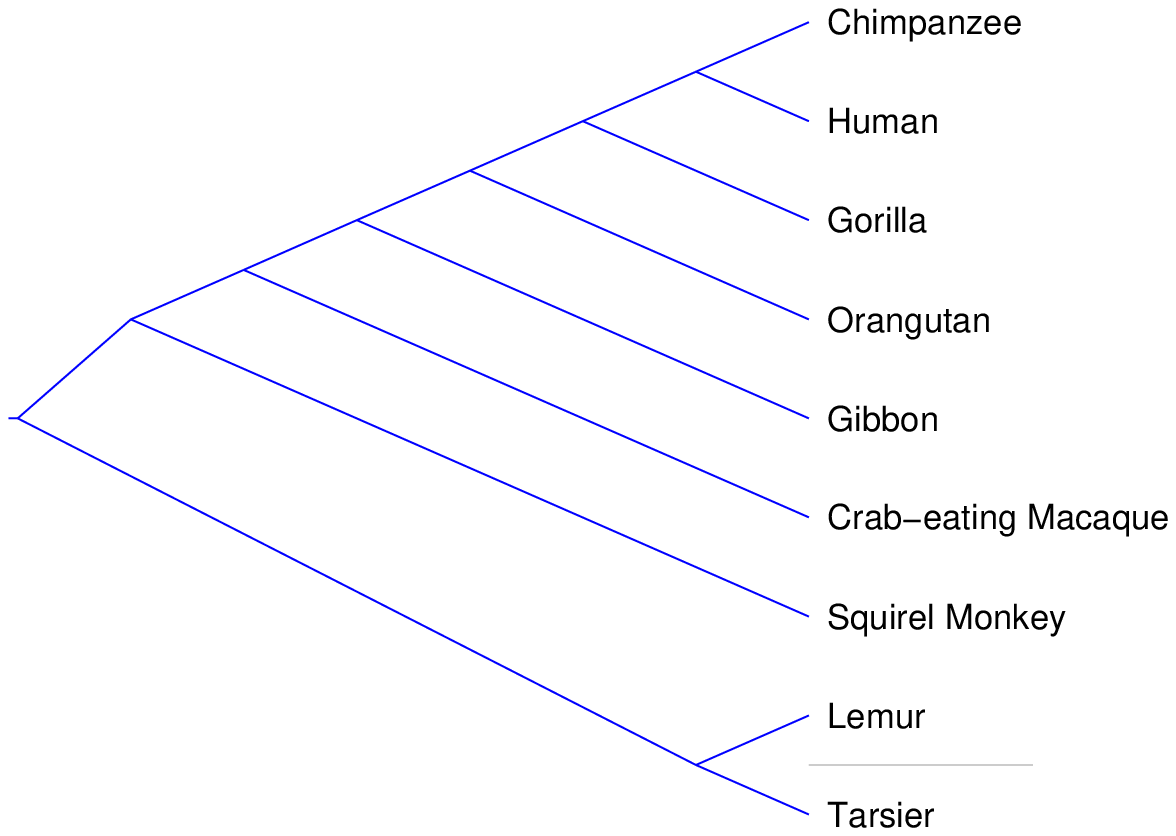}}\\
\noindent\underline{Graph 7}: The most probable phylogenetic tree topology in the primates dataset.
\end{center}

\end{document}